\documentclass[11pt]{article}
\usepackage{amstext,amssymb,amsmath,amsbsy}

\usepackage{hyperref}
\usepackage{amscd}
\usepackage{amsfonts}
\usepackage{indentfirst} 
\usepackage{verbatim}
\usepackage{amsmath}
\usepackage{amsthm} 
\usepackage{enumerate}
\usepackage{graphicx} 
\usepackage{color}
\usepackage[OT1]{fontenc}
\usepackage[latin1]{inputenc}
\usepackage[english]{babel}
\usepackage{amssymb,MnSymbol}
\usepackage{tikz}
\usepackage{dsfont}
\newtheorem{theorem}{Theorem}
\newtheorem{lemma}{Lemma}

\newcommand{\mint}{\strokedint}

\newcommand{\mR}{\mathbb{R}}
\newcommand{\mS}{\mathbb{S}}

\numberwithin{equation}{section}

\title{A refined estimate for the topological degree}

\author{
 Hoai-Minh Nguyen\footnote{EPFL SB MATHAA CAMA, Station 8,  CH-1015 Lausanne,  Switzerland, hoai-minh.nguyen@epfl.ch}
}

\begin{document}

\maketitle
%
%\vspace{-0.7cm}
%\begin{center}
%{\it To the memory of Ennio De Giorgi with emotion and admiration}
%\end{center}
%
%
\begin{abstract}
We sharpen an estimate of \cite{BBNg1} for the topological degree of continuous maps 
from a sphere $\mS^d$ into itself  in the case $d \ge 2$. This provides the answer for $d \ge 2$ to a question raised by Brezis. The problem is still open for $d=1$. 
\end{abstract}

%\end{frontmatter}
%\bigskip

%\noindent Scientific chapter:  12. Real Variable(s) Functions. 
%
%\noindent Keywords: Sobolev spaces, BV functions, non-local approximations, maximal functions. 
%
%\noindent AMS Subject classification: 46E35, 46E30, 26D15.
\medskip
{\tt AMS classification:} 47H11, 55C25, 58C35.

{\tt Keywords:} topological degree,  fractional Sobolev spaces.

\section{Introduction}

Motivated by the theory of Ginzburg Landau equations (see, e.g., \cite{BBH}), Bourgain, Brezis, and the author established in \cite{BBNg1}:
\begin{theorem} Let $d \ge 1$. For every $0 < \delta < \sqrt{2}$, there exists a positive constant $C(\delta)$ such that, for all $g \in C(\mS^d, \mS^d)$,  
\begin{equation} \label{BBN}
|\mathrm{deg}\, g| \leq C (\delta) \mathop{\int_{\mS^d} \int_{\mS^d}
}_{|g(x)-g(y)| > \delta} \frac{1}{|x-y|^{2 d}} \, dx \, dy.
\end{equation}
\end{theorem}

Here and in what follows, for $x \in \mR^{d+1}$,  $| x |$ denotes its Euclidean norm in $\mR^{d+1}$.

The constant $C(\delta)$ depends also on $d$ but for simplicity of notation we omit $d$. Estimate~\eqref{BBN} was initially suggested by Bourgain, Brezis, and Mironescu in \cite{BBM1}. It was proved in
\cite{BBM2} in the case where $d=1$ and $\delta$ is sufficiently small. 
In \cite{Ngdegree1},   the author improved \eqref{BBN} by establishing  that \eqref{BBN} holds for $0 < \delta < \ell_d = \sqrt{2 + \frac{2}{d+1}}$ with a constant $C(\delta)$ independent of $\delta$.  It was also shown there that   \eqref{BBN} does not hold for $\delta \ge \ell_d$. 

This note is concerned with the behavior of $C(\delta)$ as $\delta \to 0$. Brezis \cite{B2}  (see also \cite[Open problem 3]{B1}) conjectured that \eqref{BBN} holds with 
\begin{equation}\label{B-question}
C(\delta) = C \delta^d,  
\end{equation}
for some positive constant $C$ depending only on $d$. This conjecture is somehow motivated by the fact that \eqref{BBN}-\eqref{B-question} holds ``in the limit" as $\delta \to 0$. More precisely, it is known  that (see \cite[Theorem 2]{NgSob1})
$$
\lim_{\delta \to 0} \mathop{\int_{\mS^d} \int_{\mS^d}}_{|g(x) - g(y)| > \delta} \frac{\delta^d}{|x - y|^{2d}} \, dx \, dy = K_d \int_{\mS^d} |\nabla g (x)|^d \, dx \mbox{ for } g \in C^1(\mS^d) 
$$
for some positive constant $K_d$ depending only on $d$ and that 
$$
\mbox{deg } g  = \frac{1}{|\mS^d| }\int_{\mS^d} \mbox{Jac} (g) \mbox{ for } g \in C^1(\mS^d, \mS^d), 
$$
by Kronecker's formula. 

\medskip 
In this note, we confirm Brezis' conjecture for $d \ge 2$. The conjecture is still open for $d = 1$. Here is the result of the note.

\begin{theorem}\label{thm1} Let $d \ge 2$. There exists a positive constant $C = C(d)$,  depending only on $d$,  such that, for all $g \in C(\mS^d, \mS^d)$, 
\begin{equation}\label{thm1-conclusion2}
|\mbox{deg } g| \le C  \mathop{\int_{\mS^d} \int_{\mS^d}}_{|g(x) - g(y)| > \delta} \frac{\delta^d}{|x - y|^{2d}} \, dx \, dy  \quad \mbox{ for } 0 < \delta < 1.  
\end{equation}
\end{theorem}

\section{Proof of Theorem~\ref{thm1}}

The proof of Theorem~\ref{thm1} is in the spirit of  the approach  in \cite{BBNg1, Ngdegree1}. One of the new ingredients of the proof is the following result \cite[Theorem 1]{nguyen11}, which has its roots in \cite{BourNg}:
\begin{lemma}\label{lem1} Let $d \ge 1$, $p\ge 1$, let $B$ be  an open ball in $\mR^d$, and let $f$ be a real bounded measurable function defined in $B$. We have, for all $\delta >0$, 
\begin{equation}\label{lem1-part2}
\frac{1}{|B|^2}\int_B \int_B |f(x) - f(y)|^p \, dx \, dy \le C_{p, d} \left( |B|^{\frac{p}{d} -1}\mathop{\int_{B} \int_{B}}_{ |f(x) - f(y)| > \delta}\frac{\delta^p}{|x - y|^{d + p}} \, dx \, dy  + \delta^p \right), 
\end{equation}
for some positive constant $C_{p, d}$ depending only on $p$ and $d$. 
\end{lemma}

In Lemma~\ref{lem1}, $|B|$ denotes the Lebesgue measure of $B$. 

%\begin{remark} \rm Without imposing the condition $|x - y| < M \ell \delta$ under the integral in the RHS of \eqref{lem1-part2}, Lemma~\ref{lem1} is  given in \cite[Theorem 1]{nguyen11}. Nevertheless, the proof there also gives the assertion stated here by noting \cite[(2.4) and (2.6) in the proof of Lemma 5]{nguyen11}. 
%\end{remark}

\medskip 

We are ready to present 

\medskip 
\noindent{\bf Proof of Theorem~\ref{thm1}.}  We follow the strategy in \cite{BBNg1,Ngdegree1}. We first assume in addition that $g \in C^1(\mS^d, \mS^d)$.  Let $B$ be the open unit ball in $\mR^{d+1}$ and  
let $u: B \to B$ be the average extension of $g$, i.e., 
\begin{equation}
u(X)   =   \mint_{B(x, r)} g(s) \, ds \mbox{ for } X \in B, 
\end{equation}
where $x = X/ |X|$,  $r = 2 (1 - |X|)$, and $B(x, r): = \big\{y \in \mS^d; |y-x| \le r \big\}$. In this proof, $\mint_{D} g(s) \, ds $ denotes the equantity $ \frac{1}{|D|} \int_{D} g(s) \, ds$ for a measurable subset $D$ of $\mS^{d}$ with positive ($d$-dimensional Hausdorff) measure. 
Fix $\alpha  = 1/ 2$  and for every $x \in \mS^d$, let $\rho(x)$ be the length of the largest radial interval coming from $x$ on which $|u| > \alpha$ (possibly $\rho(x)$ = 1). In particular,  if $\rho(x) < 1$,  then
\begin{equation}\label{rho}
\left| \mint_{B(x, 2\rho(x))} g(s) \, ds \right| = 1/2. 
\end{equation}
By \cite[(7)]{BBNg1}, we have
\begin{equation}\label{thm1-part1}
|\mbox{deg } g| \le C\mathop{\int_{\mS^d}}_{\rho(x) < 1} \frac{1}{\rho(x)^d} \, dx.
\end{equation}
Here and in what follows, $C$ denotes a positive constant which is independent of $x,  \,  \xi \, , \eta$,  $g$, and $\delta$, and can change from one place to another. 

\medskip
We now implement ideas involving Lemma~\ref{lem1} applied with $p=1$. 
We have, by \eqref{rho}, 
$$
\mint_{B(x, 2\rho(x))} \mint_{B(x, 2\rho(x))} |g(\xi) - g(\eta)| \, d \xi \,  d \eta \ge \mint_{B(x, 2\rho(x))} \left|g(\xi) - \mint_{B(x, 2\rho(x))} g(\eta) \, d \eta \right| \, d \xi  \ge C. 
$$
This yields, for some $1 \le j_0 \le d+1$, 
$$
\mint_{B(x,2 \rho(x))} \mint_{B(x, 2\rho(x))} |g_{j_0} (\xi) - g_{j_0} (\eta)| \, d \xi \,  d \eta \ge  C, 
$$
where $g_j$ denotes the $j$-th component of $g$.  It follows from \eqref{lem1-part2} that, for some $\delta_0 > 0$ ($\delta_0$ depends only on $d$)  and for $0< \delta < \delta_0$, 
\begin{equation*}
 \rho(x)^{1-d} \mathop{\int_{B(x,2 \rho(x))} \int_{B(x, 2\rho(x))}}_{|g_{j_0}(\xi) - g_{j_0}(\eta)| >  \delta} \frac{\delta}{|\xi-\eta|^{d+ 1}} \, d \xi \, d \eta  \ge C, 
\end{equation*}
which implies 
\begin{equation}\label{lem1-part2-*}
\sum_{j=1}^{d+1} \rho(x)^{1-d} \mathop{\int_{B(x,2 \rho(x))} \int_{B(x, 2\rho(x))}}_{ |g_j(\xi) - g_j(\eta)| >  \delta} \frac{\delta}{|\xi-\eta|^{d+ 1}} \, d \xi \, d \eta  \ge C. 
\end{equation}
 Since 
$$
 \rho(x)^{1-d} \mathop{\int_{B(x, 2 \rho(x))} \int_{B(x, 2\rho(x))}}_{|\xi - \eta| > C_1 \rho(x) \delta} \frac{\delta}{|\xi-\eta|^{d+ 1}} \, d \xi \, d \eta < \frac{C}{2(d+1)}, 
$$
if $C_1>0$ is large enough (the largeness of $C_1$ depends only on $C$ and $d$), it follows from~\eqref{lem1-part2-*} that 
\begin{equation}\label{thm1-part2}
\sum_{j=1}^{d+1} \rho(x)^{1-d} \mathop{\mathop{\int_{B(x, 2\rho(x))} \int_{B(x, 2\rho(x))}}_{ |g_j(\xi) - g_j(\eta)| >  \delta}}_{|\xi - \eta| \le C \rho(x)\delta } \frac{\delta}{|\xi-\eta|^{d+ 1}} \, d \xi \, d \eta  \ge C.  
\end{equation}
We derive from \eqref{thm1-part1} and \eqref{thm1-part2} that, for $0< \delta < \delta_0$, 
\begin{equation*}
|\mbox{deg } g| \le C \mathop{\int_{\mS^d}}_{\rho(x) < 1} \frac{1}{\rho(x)^{2d - 1}} \, dx \sum_{j=1}^{d+1} \mathop{\mathop{\int_{B(x, 2\rho(x))} \int_{B(x, 2\rho(x))}}_{ |g_j(\xi) - g_j(\eta)| >  \delta}}_{|\xi - \eta| \le C \rho(x) \delta} \frac{\delta}{|\xi-\eta|^{d+ 1}} \, d \xi \, d \eta. 
\end{equation*}
This implies, by Fubini's theorem,  that,  for $0< \delta < \delta_0$,  
\begin{equation}\label{thm1-part3}
|\mbox{deg } g| \le   C  \sum_{j=1}^{d+1} \mathop{\int_{\mS^d} \int_{\mS^d}}_{|g_j(\xi) - g_j(\eta)| >  \delta} \frac{\delta}{|\xi-\eta|^{d+ 1}} \, d \xi \, d \eta \mathop{\mathop{\int}_{\rho(x) \ge C|\xi - \eta|/ \delta}}_{2\rho(x) > |x - \xi|} \frac{1}{\rho(x)^{2d - 1}} \, dx.
\end{equation}
We  have 
\begin{align*}
\mathop{\mathop{\int}_{2 \rho(x) > |x - \xi|}}_{\rho(x) \ge C |\xi - \eta|/ \delta} \frac{1}{\rho(x)^{2d - 1}} \, dx & \le   \mathop{\mathop{\int}_{2\rho(x) > |x - \xi|}}_{ |x - \xi|> C |\xi - \eta|/ \delta} \frac{1}{\rho(x)^{2d - 1}} \, dx + \mathop{\mathop{\int}_{  \rho(x) \ge C |\xi - \eta|/ \delta }}_{|x - \xi|  \le C  |\xi - \eta|/\delta } \frac{1}{\rho(x)^{2d - 1}} \, dx \\[6pt]
& \le   \mathop{\int}_{ |x - \xi|> C |\xi - \eta|/ \delta} \frac{C}{|x - \xi|^{2d - 1}} \, dx + \mathop{\int}_{|x - \xi|  \le C  |\xi - \eta|/\delta } \frac{C \delta^{2d-1}}{|\xi - \eta|^{2d - 1}} \, dx. 
\end{align*}
Finally, we use the assumption that $d \ge 2$. Since $d > 1$, it follows that  
\begin{equation}\label{thm1-part4}
\mathop{\mathop{\int}_{\rho(x) > |x - \xi|}}_{\rho(x) \ge C |\xi - \eta|/ \delta} \frac{1}{\rho(x)^{2d - 1}} \, dx \le \frac{C \delta^{d-1} }{ |\xi - \eta|^{d-1}}.  
\end{equation}
Combining \eqref{thm1-part3} and \eqref{thm1-part4} yields, for $0 < \delta < \delta_0$, 
\begin{equation}\label{thm1-part5}
|\mbox{deg } g| \le  C  \sum_{j=1}^{d+1}  \mathop{\int_{\mS^d} \int_{\mS^d}}_{ |g_j(\xi) - g_j(\eta)| >  \delta} \frac{\delta^d}{|\xi-\eta|^{2d}} \, d \xi \, d \eta. 
\end{equation}
Assertion~\eqref{thm1-conclusion2} is now a direct consequence of  \eqref{thm1-part5} for $\delta < \delta_0$ and \eqref{BBN} for $\delta_0 \le \delta < 1$. 

\medskip 
The proof in the case $g \in C(\mS^d, \mS^d)$ can be derived from the case $g \in C^1(\mS^d, \mS^d)$ via a standard approximation argument. The details are omitted. 
\qed

\bigskip 
\noindent {\bf Acknowledgement:} The author warmly  thanks Haim Brezis  for  communicating \cite{B2} and Haim Brezis and Itai Shafrir  for interesting discussions. 

%\begin{remark} \rm Similar estimates as in \eqref{lem1-part2-*} for $\delta$ of the order 1 and $p=d$ are used in \cite{BBNg1} and \cite{Ngdegree1}. 
%Nevertheless, this kind of estimates do not fit in the setting of Theorem~\ref{thm1}. 
%\end{remark}

%\begin{remark} \rm Let $I$ be a bounded interval of $\mR$ and let $f$ be a real measurable function defined in $I$. 
%It is known that \cite[Theorem 1]{nguyen11}
%$$
%\mint_I \mint_I |f(x) - f(y)| \, dx \, dy \le C \left(\mathop{\int_{I} \int_{I}}_{ |f(x) - f(y)| > \delta} \frac{1}{|x - y|^{2}} \, dx \, dy  + \delta \right), 
%$$
%for some universal positive constant $C$. Nevertheless, we are not able to use this estimate to obtain \eqref{thm1-conclusion2} in the one dimensional case. 
%
%\end{remark}

\end{document}